\documentclass[final]{siamltex}

\usepackage{amsmath, amssymb, amsfonts}
\usepackage{graphicx}

\newcommand{\eq}[1]{(\ref{#1})}
\newcommand{\eps}{\epsilon}

\newcommand{\Rl}{\mathbb{R}}

\newcommand{\moddf}{|\partial_x|}

\newcommand{\df}{\partial_x}

\newcommand{\Hi}{\operatorname{\mathbf{H}}}
\newcommand{\Id}{\operatorname{\mathbf{I}}}

\title{Enhanced Lifespan of Smooth Solutions of a Burgers-Hilbert Equation\thanks{Submitted: September 29, 2011.}}

\author{John K. Hunter \thanks{Department of Mathematics,University of California at Davis.
        Partially supported by the NSF under grant number DMS-0072343.}
        \and
        Mihaela Ifrim \thanks{Department of Mathematics, University of California at Davis}.}
        
\date{September 29, 2011}

\begin{document}

\maketitle

\begin{abstract}
We consider an initial value problem for a quadratically nonlinear inviscid
Burgers-Hilbert equation that models the motion of vorticity discontinuities. We use a
normal form transformation, which is implemented by means of a near-identity
coordinate change of the independent spatial variable, to prove the existence of small, smooth solutions
over cubically nonlinear time-scales. For vorticity discontinuities, this result means that
there is a cubically nonlinear time-scale before the onset of filamentation.
\end{abstract}

\begin{keywords}
Normal form transformations, nonlinear waves, inviscid Burgers equation, vorticity discontinuities.
\end{keywords}

\begin{AMS}
37L65, 76B47.
\end{AMS}

\pagestyle{myheadings}
\thispagestyle{plain}
\markboth{HUNTER AND IFRIM}{BURGERS-HILBERT EQUATION}

\section{Introduction}
We consider the following initial value problem for an inviscid Burgers-Hilbert equation for $u(t,x; \eps)$:
\begin{align}
\label{bheq}
\begin{split}
&u_t+ \eps u u_x = \mathbf{H}\left[ u\right] ,\\
&u(0,x;\eps) = u_0(x).
\end{split}
\end{align}
In \eq{bheq}, $\mathbf{H}$ is the spatial Hilbert transform, $\eps$ is a small parameter, and
$u_0$ is given smooth initial data. This Burgers-Hilbert equation is a model equation for
nonlinear waves with constant frequency \cite{biello}, and it provides an effective equation for the motion of a vorticity discontinuity
in a two-dimensional flow of an inviscid, incompressible fluid \cite{biello, marsden}.
Moreover, as shown in \cite{biello}, even though \eq{bheq} is quadratically nonlinear it provides a formal
asymptotic approximation for the small-amplitude motion of a planar vorticity discontinuity
located at $y = \eps u(t,x;\eps)$ over cubically nonlinear
time-scales.

We assume for simplicity that $x\in \mathbb{R}$, in which case
the Hilbert transform is given by
\[
\mathbf{H}[u](t,x;\eps) = \mathrm{p.v.} \frac{1}{\pi} \int\frac{u(t,y;\eps)}{x-y}\, dy.
\]
We will show that smooth solutions of \eq{bheq} exist for times of the order $\eps^{-2}$ as $\eps \to 0$. Explicitly,
if $H^s(\Rl)$ denotes the standard Sobolev space of functions with $s$ weak $L^2$-derivatives, we prove the following result:

\begin{theorem}
\label{th:main}
Suppose that $u_0 \in H^2(\Rl)$.
There are constants $k >0$ and $\eps_{0} > 0$, depending only on $\Vert u_{0}\Vert _{H^2}$,
such that for every $\eps$ with
$|\eps| \leq \eps _{0}$ there exists a solution
\[
u \in C\left( I^\eps;{H}^{2}\left( \mathbb{R}\right) \right) \cap C^1 \left(I^\eps;{H}^{1}\left(\mathbb{R}\right) \right)
\]
of (\ref{bheq})
defined on the time-interval $I^\eps=\left[-{k}/{\eps^2},{k}/{\eps^2}\right]$.
\end{theorem}

The cubically nonlinear $O(\eps^{-2})$ lifespan of smooth solutions for the Burgers-Hilbert equation is
longer than the quadratically nonlinear $O(\eps^{-1})$ lifespan for the inviscid Burgers equation
$u_t + \eps uu_x=0$.
The explanation of this enhanced lifespan is that the quadratically nonlinear term of the order $\eps$ in \eq{bheq}
is nonresonant for the linearized equation. To see this, note that
the solution of the linearized equation $u_t = \mathbf{H}[u]$ is given by $u = e^{t\mathbf{H}} u_0$, or
\[
u(t,x) = u_0(x) \cos t + h_0(x) \sin t,\qquad h_0 = \mathbf{H}[u_0],
\]
as may be verified by use of the identity $\Hi^2 = -\Id$.
This solution oscillates with frequency one between the initial data and its
Hilbert transform, and the effect of the nonlinear forcing term $\eps u u_x$ on the linearized equation
averages to zero because it contains no Fourier component in time whose frequency is equal to one.

\begin{figure}[htp]
\centering
\includegraphics[width=4in]{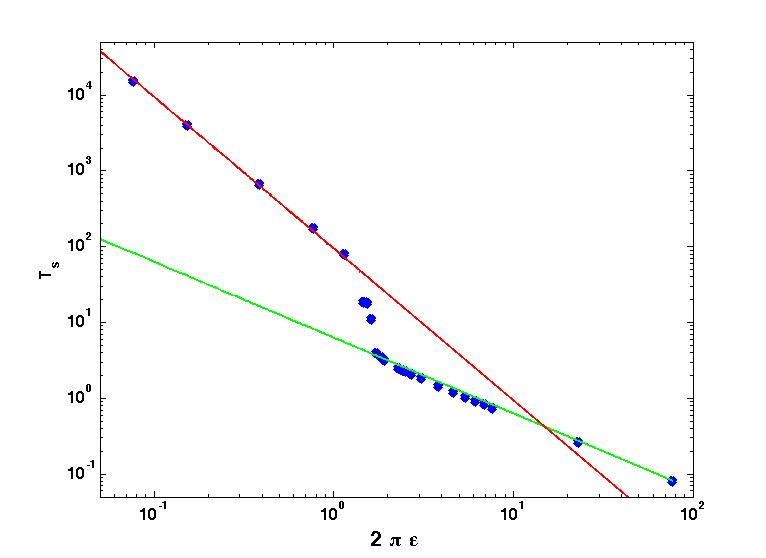}
\caption{Logarithm of the singularity formation time $T_s$ for the Burgers-Hilbert equation \eq{bheq} versus the logarithm
of $2 \pi \eps$ for fixed initial data $u_0$. Numerical solutions are shown by diamonds. The steeper line is a
formal asymptotic prediction from \cite{biello} for $\eps\ll 1$, which gives $T_s = 2.37\, \eps^{-2}$. The shallower line is the
singularity formation time for the  inviscid Burgers equation, which gives $T_s = \eps^{-1}$. (See \cite{biello} for further details.)}
\label{u_sing}
\end{figure}

Alternatively, one can view the averaging of the nonlinearity as a consequence of the fact that the nonlinear steepening
of the profile in one phase of the oscillation is canceled by its expansion in the other phase. This phenomenon is
illustrated by numerical results from \cite{biello}, which are reproduced in Figure~\ref{u_sing}. The transition from
an $O(\eps^{-1})$ lifespan for large $\eps$ to an $O(\eps^{-2})$ lifespan for small $\eps$ is remarkably rapid: once a
singularity fails to form over the first oscillation in time, a smooth solution typically persists over many oscillations.

In the context of the motion of a vorticity discontinuity,
%the boundary of a vortex patch \cite{chemin,majda} may be locally approximated by such a discontinuity.
the formation of a singularity in a solution of \eq{bheq} corresponds to the filamentation
of the discontinuity \cite{dritschel}. The result proved here corresponds to an enhanced lifespan before
nonlinear `breaking' of the discontinuity leads to the formation of a filament.

There are three main difficulties in the proof of Theorem~\ref{th:main}. The first is that the presence
of a quadratically nonlinear term in \eq{bheq} means that straightforward energy estimates prove the existence of smooth solutions
only on time-scales of the order $\epsilon^{-1}$. Following the idea introduced by Shatah \cite{shatah} in the context of PDEs,
and used subsequently by other authors,
we remove the quadratically nonlinear term of the order $\eps$ by a normal form or near-identity transformation,
replacing it by a cubically nonlinear term of the order $\eps^2$. The second difficulty is that a
standard normal form transformation of the dependent variable,
of the type used by Shatah, leads to a loss of spatial derivatives because we are using a lower-order linear term $\mathbf{H}[u]$ to eliminate
a higher-order nonlinear term $\eps uu_x$. The third difficulty is that \eq{bheq} is nondispersive and solutions of the linearized equation oscillate
but do not decay in time. Thus, we cannot use any dispersive decay in time to control the loss of spatial derivatives.

The key idea in this paper that avoids these difficulties is to make a transformation of the independent variable, rather than the dependent variable.
We write
\begin{equation}
h(t,x;\eps) = \mathbf{H}[u](t,x;\eps)
\label{defh}
\end{equation}
and define
\begin{equation}
g(t,\xi;\eps) = h(t,x;\eps),\qquad x = \xi - \eps g(t,\xi;\eps).
\label{nearid_trans}
\end{equation}
Then, as we will show, the transformed function $g(t,\xi;\eps)$ satisfies an integro-differential equation of the form
\begin{align}
\begin{split}
g_{t}(t,\xi;\eps) &= \mathrm{p.v.} \frac{1}{\pi} \int\frac{g(t,\tilde{\xi};\eps)}{\xi-\tilde{\xi}}\, d\tilde{\xi}
\\
&\qquad - \frac{1}{\pi} \eps^2 \partial_{\xi} \int  (\xi-\tilde{\xi}) g_{\tilde{\xi}}(t,\tilde{\xi};\eps)
\phi\left(\frac{g(t,\xi;\eps)-g(t,\tilde{\xi};\eps)}{\xi-\tilde{\xi}}; \eps\right)\, d\tilde{\xi}
\end{split}
\label{g-eq}
\end{align}
where $\phi(c;\eps)$ is a smooth function, given in Lemma~\ref{lemma:geq}.
The term of the order $\eps$ has been removed from \eq{g-eq}, and the equation has good energy estimates
that imply the enhanced lifespan of smooth solutions.

The interpretation of the transformation \eq{nearid_trans} is not entirely clear. On taking the Hilbert transform of \eq{bheq}
we get $h_t = - u + O(\eps)$, so that
\[
x_t = -\eps g_t = -\eps h_t + O(\eps^2) = \eps u + O(\eps^2).
\]
Thus the transformation $\xi\mapsto x$ in \eq{nearid_trans}
agrees up to the order $\eps$ with a transformation from characteristic to
spatial coordinates for \eq{bheq}. The coordinate $\xi$, however, differs from $x$ even when $t=0$, and
the use of characteristic coordinates does not appear to simplify the analysis.
As a partial motivation for \eq{nearid_trans}, we show in Section~\ref{sec:normal_form} that
it agrees to leading
order in $\eps$ with a normal form transformation of the dependent variable that is given in \cite{biello}.
We were not able, however, to use the latter normal form transformation
to prove Theorem~\ref{th:main} because of the loss of derivatives in the higher-order terms.
%the transformation \eq{nearid_trans} avoids this difficulty.

We consider \eq{bheq} on the real line for simplicity. Equation \eq{bheq} is nondispersive and solutions of the linearized
equation oscillate in time. Thus, our proof does not depend on any dispersive decay of the solutions in time,
and a similar result would apply to spatially periodic solutions. Theorem~\ref{th:main} is also presumably
true in $H^s$ for any $s > 3/2$; we consider $s=2$ to avoid complications associated with the use of fractional derivatives.
A proof of singularity formation for \eq{bheq} under certain conditions on $u_0$ and $\eps$
is given in \cite{castro}.

\section{Proof of the Theorem}
\label{sec:proof}

In this section, we prove Theorem~\ref{th:main}.
It follows from standard energy arguments (\textit{e.g.}\ \cite{kato}) that \eq{bheq} has a unique local $H^2$-solution in a time-interval $J^\eps$
depending on the $H^2$-norm of the initial data and $\eps$. Moreover, for any $s\ge 2$, the solution remains in $H^s$ if the initial data is in $H^s$
and depends continuously on the initial data in $C(J^\eps; H^s)$. Thus, in order to prove Theorem~\ref{th:main}
it is sufficient to prove an \textit{a priori} $H^2$-bound for smooth solutions $u \in C^\infty\left(I^\eps; H^\infty(\Rl)\right)$
where $H^\infty(\Rl) = \cap_{s=1}^\infty H^s(\Rl)$. To derive this bound, we first transform the equation to remove the order $\eps$ term
and then carry out $H^2$-estimates on the transformed equation. The required computations, such as integrations by parts, are
justified for these smooth solutions.

\subsection{Near-identity transformation}

Let $h$ denote the Hilbert transform of $u$, as in \eq{defh}. Taking the Hilbert transform of \eq{bheq}, using the identity
\[
\mathbf{H}\left[u^2 - h^2\right] = 2 h u
\]
and the fact that $u= - \mathbf{H}[h]$, we find that $h$ satisfies the equation
\begin{equation}
\label{htrans}
h_{t}+\eps \left\lbrace \mathbf{H}\left[hh_{x}\right]-h\mathbf{H}\left[ h_{x}\right] -\mathbf{H}\left[ h\right]h_{x}\right\rbrace =\mathbf{H}\left[ h\right].
\end{equation}
We will make the change of variables \eq{nearid_trans} in \eq{htrans}, so first we discuss \eq{nearid_trans}.

The map $\xi\mapsto x$ is smoothly invertible if
$|\eps g_\xi| < 1$, which holds by Sobolev embedding if $\|\eps g\|_{H^2}$ is sufficiently small.
Specifically, we have  the Gagliardo-Nirenberg-Moser inequality
\begin{equation}
\|g_\xi\|_{L^\infty} \le N \|g\|^{1/4}_{L^2} \|g_{\xi\xi}\|^{3/4}_{L^2},
\label{gag_nir}
\end{equation}
where we can take, for example,
\[
N = \sqrt{\frac{8}{3}}.
\]
We assume throughout this section that
\begin{equation}
N\|\eps g\|^{1/4}_{L^2} \|\eps g_{\xi\xi}\|^{3/4}_{L^2} \le \frac{1}{2},
\label{c2small}
\end{equation}
which ensures that
\begin{equation}
\|\eps g_\xi\|_{L^\infty} \le \frac{1}{2}.
\label{cinftysmall}
\end{equation}

By the chain rule,
\begin{equation}
h_{x}=\frac{g_{\xi}}{1-\eps g_{\xi}},
\qquad
h_{xx}=\frac{g_{\xi\xi}}{(1-\eps g_{\xi})^3}.
\label{hx}
\end{equation}
Thus, if \eq{cinftysmall} holds, then
\[
\int_{\Rl} h^2 \, dx = \int_{\Rl} g^2 \left(1-\eps g_\xi\right)\, d\xi
= \int_{\Rl} g^2 \, d\xi,
\qquad
\int_{\Rl} h_{xx}^2 \, dx = \int_{\Rl} \frac{g_{\xi\xi}^2}{\left(1-\eps g_\xi\right)^5}\, d\xi.
\]
Hence, since $\mathbf{H}$ is an isometry on $H^s$,
\begin{equation}
\|u\|_{L^2} = \|g\|_{L^2},
\qquad
\left(\frac{2}{3}\right)^{5/2} \|g_{\xi\xi}\|_{L^2}
\le
\|u_{xx}\|_{L^2}
\le
2^{5/2}  \|g_{\xi\xi}\|_{L^2},
\label{normest}
\end{equation}
and $H^2$-estimates for $g$ imply $H^2$-estimates for $u$.

Conversely, one can use the contraction mapping theorem on $C_0(\Rl)$, the space
of continuous functions that decay to zero at infinity, to show that if $h_0\in C_0^1(\Rl)$ and
\begin{equation}
\|\eps h_{0x}\|_{L^\infty} < 1
\label{hginvert}
\end{equation}
then there exists a function $g_0(\cdot;\eps)\in C_0(\Rl)$
such that
\[
h_0\left(\xi - \eps g_0(\xi;\eps)\right) = g_0(\xi;\eps).
\]
The function $g_0$ is smooth if $h_0$ is smooth, and $\|\eps g_{0\xi}\|_{L^\infty} \le 1/2$ if
$\|\eps h_{0x}\|_{L^\infty} \le 1/3$.
Thus, we can obtain initial data for $g$ from the initial data for $h$.

From \eq{nearid_trans}, we have
\begin{equation*}
h_{t}=\frac{g_{t}}{1-\eps g_{\xi}},
\qquad
\mathbf{H}\left[h\right]= \mathrm{p.v.}\frac{1}{\pi}\int_{\Rl}\left[ \frac{1-\eps\tilde{g}_{\tilde{\xi}}}{\xi-\tilde{\xi}-\eps(g-\tilde{g})}\right]\tilde{g} \, d\tilde{\xi}
\end{equation*}
where we use the notation
\[
g=g(t,\xi;\eps),\qquad \tilde{g}=g(t,\tilde{\xi};\eps).
\]
Using these expressions, together with \eq{hx}, in (\ref{htrans}) and simplifying the result, we find that $g(t,\xi;\eps)$ satisfies the
following nonlinear integro-differential equation:
\[
g_{t}= \mathrm{p.v.}\frac{1}{\pi}\int_{\Rl}\frac{\tilde{g}+\eps(g-2\tilde{g})\tilde{g}_{\tilde{\xi}}
-\eps^2 (g-\tilde{g})g_{\xi}\tilde{g}_{\tilde{\xi}}}{\xi-\tilde{\xi}-\eps (g-\tilde{g})}\, d\tilde{\xi}.
\]
Subtracting off the leading order term in $\eps$ from the integrand, we may write this equation as
\begin{equation}
\label{f}
g_{t}=\mathbf{H}[g]
+ \frac{1}{\pi}\eps^2 \int_{\Rl}\left(\frac{g-\tilde{g}}{x-\tilde{x}}\right)\left\{
\left(\frac{\tilde{g}}{\xi-\tilde{\xi}}\right)\left[\frac{g-\tilde{g}}{\xi-\tilde{\xi}} - \tilde{g}_{\tilde{\xi}}\right]
+ \tilde{g}_{\tilde{\xi}}\left[\frac{g-\tilde{g}}{\xi-\tilde{\xi}} - g_\xi\right]\right\}\, d\tilde{\xi}
\end{equation}
where
\[
\mathbf{H}[g](t,\xi;\eps) = \mathrm{p.v.}\frac{1}{\pi} \int_{\Rl} \frac{g(t,\tilde{\xi};\eps)}{\xi - \tilde{\xi}} \, d\tilde{\xi}
\]
denotes the Hilbert transform of $g$ with respect to $\xi$ and
\[
x = \xi - \eps g(t,\xi;\eps),\qquad \tilde{x} = \tilde{\xi} - \eps g(t,\tilde{\xi};\eps).
\]
The integral of the order $\eps^2$ in \eq{f} is not a principal value integral
since the integrand is a smooth function of $(\xi,\tilde{\xi})$.
Finally, we observe that this equation can be put in the form \eq{g-eq}.

\begin{lemma}
\label{lemma:geq}
An equivalent form of equation (\ref{f}) is given by
\begin{equation}
\label{finalfinaleq}
g_{t}=\mathbf{H}\left[g \right] -\frac{1}{\pi}\eps ^2 \partial_{\xi}\int_{\Rl}(\xi-\tilde{\xi})\tilde{g}_{\tilde{\xi}}\,
\phi \left( \frac{g-\tilde{g}}{\xi-\tilde{\xi}};\eps\right) \, d\tilde{\xi},
\end{equation}
where
\begin{equation}
\phi(c;\eps)= - \frac{1}{\eps^2}\left\lbrace \log \left( 1-\eps c\right) + \eps c \right\rbrace.
\label{defphi}
\end{equation}
\end{lemma}

\begin{proof}
First, we check that \eq{finalfinaleq} is well-defined.
Abusing notation slightly, we write
\begin{equation}
c =  \frac{g-\tilde{g}}{\xi-\tilde{\xi}}.
\label{defc}
\end{equation}
From \eq{defphi},
\begin{equation}
\phi_c(c;\eps) = \frac{c}{1-\eps c}, %\qquad \phi_c\left(\frac{g-\tilde{g}}{\xi-\tilde{\xi}};\eps\right) = \frac{g-\tilde{g}}{x-\tilde{x}}.
\label{defphic}
\end{equation}
so $|\phi(c;\eps)| \le c^2$ when $|\eps c| \le 1/2$, which is implied by \eq{cinftysmall}.
In that case
\[
\left|\int_{\Rl}(\xi-\tilde{\xi})\tilde{g}_{\tilde{\xi}}\,
\phi \left(c;\eps\right) \, d\tilde{\xi}\right| \le \int_{\Rl} \left|\left(g-\tilde{g}\right) \tilde{g}_{\tilde{\xi}} c\right|\,d\tilde{\xi}.
\]
We use $|g-\tilde{g}| \le 2 \|g\|_{L^\infty}$ in the right hand side of this inequality and
apply the Cauchy-Schwartz inequality to get
\begin{equation}
\left|\int_{\mathbb{R}}(\xi-\tilde{\xi})\tilde{g}_{\tilde{\xi}}\,
\phi \left(c;\eps\right) \, d\tilde{\xi}\right| \le 2 \|g\|_{L^\infty}\|{g}_{\xi}\|_{L^2}\left\|c\right\|_{L^2}
\label{tempint}
\end{equation}
where
\[
\left\|c\right\|_{L^2} = \left[\int_{\mathbb{R}} \left(\frac{g-\tilde{g}}{\xi-\tilde{\xi}}\right)^2\, d\tilde{\xi}\right]^{1/2}
\]
denotes the $L^2$-norm of $c$ with respect to $\tilde{\xi}$, which is a function of $\xi$.
Temporarily suppressing the $(t;\eps)$-variables
and denoting the derivative of $g$ with respect to $\xi$ by $g^\prime(\xi) = g_\xi(\xi)$, we have from
the Taylor integral formula that
\[
c = \int_0^1 g^\prime\left(\tilde{\xi} + r(\xi - \tilde{\xi})\right)\, dr,
\]
and the Cauchy-Schwartz inequality implies that
\begin{align*}
\left\|c\right\|^2_{L^2}(\xi) &= \int_{\Rl} c^2 \, d\tilde{\xi}
\\
&= \int_0^1 \int_0^1 \int_{\Rl} g^\prime\left(\tilde{\xi} + r(\xi - \tilde{\xi})\right)g^\prime\left(\tilde{\xi} + s(\xi - \tilde{\xi})\right) \, d\tilde{\xi} dr ds
\\
& \le \int_0^1 \int_0^1 \left(\int_{\Rl} g^{\prime2}\left(\tilde{\xi} + r(\xi - \tilde{\xi})\right)\, d\tilde{\xi}\right)^{1/2}
\left(\int_{\Rl} g^{\prime2}\left(\tilde{\xi} + s(\xi - \tilde{\xi})\right) \, d\tilde{\xi}\right)^{1/2} dr ds
\\
& \le \left(\int_0^1 \int_0^1 \frac{1}{\sqrt{rs}} dr ds\right) \left(\int_{\Rl} g^{\prime2}(\xi)\, d\xi\right)
\\
&\le 4 \left(\int_{\Rl} g^{\prime2}(\xi)\, d\xi\right).
\end{align*}
Thus,
\begin{equation}
\sup_{\xi\in \Rl}\left(\int_{\Rl} c^2 \, d\tilde{\xi} \right)^{1/2} \le 2 \left\|g_\xi\right\|_{L^2}.
\label{cl2}
\end{equation}
Using this estimate in \eq{tempint}, we get
\[
\sup_{\xi\in \Rl}\left|\int_{\Rl}(\xi-\tilde{\xi})\tilde{g}_{\tilde{\xi}}\,
\phi \left(c;\eps\right) \, d\tilde{\xi}\right| \le 4 \|g\|_{L^\infty} \|{g}_{\xi}\|^2_{L^2}.
\]
Thus, the $\tilde{\xi}$-integral in \eq{finalfinaleq} converges when $g \in H^1(\Rl)$ and is, in fact, a uniformly bounded function of $\xi$.

To verify that \eq{finalfinaleq} agrees with \eq{f}, we
take the $\xi$-derivative under the integral in \eq{finalfinaleq}, use \eq{defphic} which implies that
\[
\phi_c\left(c;\eps\right) = \frac{g-\tilde{g}}{x-\tilde{x}},
\]
and integrate by parts in the result. This gives
\begin{equation}
g_t = \mathbf{H}\left[g \right] + \frac{1}{\pi}\eps^2 \int_{\Rl} \left(\frac{g-\tilde{g}}{x-\tilde{x}}\right)
\left[\tilde{g} c_{\tilde{\xi}} - (\xi-\tilde{\xi}) \tilde{g}_{\tilde{\xi}} c_\xi \right]
\, d\tilde{\xi}.
\label{tempgeq}
\end{equation}
Using the equations
\begin{equation}
c_{\tilde{\xi}} =  \frac{c - \tilde{g}_{\tilde{\xi}}}{\xi-\tilde{\xi}},
\qquad
c_\xi = - \frac{c - g_{\xi}}{\xi-\tilde{\xi}} ,\qquad
\label{defcxi}
\end{equation}
in \eq{tempgeq} and comparing the result with \eq{f} proves the lemma.
\end{proof}

%------------------------------------------------------------------------------------------------------------------
\subsection{Energy Estimates}

Multiplying (\ref{finalfinaleq}) by $g$, integrating the result with respect to $\xi$, and integrating by parts with respect
to $\xi$, we find that the right-hand side
vanishes by skew-symmetry in $(\xi,\tilde{\xi})$ so that
\begin{equation}
\frac{d}{dt} \left\|g\right\|_{L^2} = 0.\label{conl2norm}
\end{equation}
The conservation of $\|g\|_{L^2}$
is consistent with the conservation
of $\|u\|_{L^2}$ from \eq{bheq}.
%The initial condition for $g(t,\xi;\eps)$ is
%\[
%g(0,\xi;\eps) = h_0\left(x(0,\xi;\eps)\right).
%\]
Hence, from \eq{normest}, we have $\|g\|_{L^2}= \|g_0\|_{L^2} = \|u_0\|_{L^2}$.

Differentiating (\ref{finalfinaleq}) twice with respect to $\xi$, multiplying the result by $g_{\xi \xi}$, integrating with respect to $\xi$, and
integrating by parts with respect to $\xi$, we get
\begin{equation}
\frac{d}{dt}\int_{\Rl}g_{\xi \xi}^2\, d\xi=\eps^2 I
\label{H2eq}
\end{equation}
where
\begin{equation}
I = \int_{\Rl^2} g_{\xi \xi \xi} \partial_{\xi}^2\left[ (\xi -\tilde{\xi})\tilde{g}_{\tilde{\xi}}\phi (c;\eps) \right] \, d\xi d\tilde{\xi}.
\label{tempI}
\end{equation}
The following lemma estimates $I$ in terms of the $H^2$-norm of $g$.

\begin{lemma}
\label{lemma:est}
Suppose that $I$ is given by \eq{tempI} where $\phi$ is defined in \eq{defphi},  and $c$ is defined in \eq{defc}.
There exists a numerical constant $A>0$ such that
\begin{equation}
|I| \le A \left\|g_\xi\right\|_{L^2} \|g_{\xi\xi}\|^3_{L^2}
\label{Iest}
\end{equation}
whenever $g\in H^\infty(\Rl)$ satisfies \eq{c2small}.
\end{lemma}

\begin{proof}
We first convert the $\tilde{\xi}$-derivative in the expression \eq{tempI} for $I$ to a $\xi$-derivative. Let
\[
\Phi^\prime(c;\eps) = \phi(c;\eps),
\]
where a prime on $\Phi$ and related functions denotes a derivative with respect to $c$.
It follows from \eq{defcxi} that
\begin{align*}
(\xi -\tilde{\xi})\tilde{g}_{\tilde{\xi}}\phi(c;\eps) &= (\xi -\tilde{\xi})\left[c - (\xi-\tilde{\xi}) c_{\tilde{\xi}}\right]\phi(c;\eps)
\\
&= (\xi -\tilde{\xi}) c\phi(c;\eps) - (\xi-\tilde{\xi})^2 \Phi_{\tilde{\xi}}(c;\eps).
\end{align*}
We use this equation in \eq{tempI} and integrate by parts with respect to $\tilde{\xi}$ in the term involving $\Phi$.
Since $g$ is independent of $\tilde{\xi}$, this gives
\begin{equation}
I =\int_{\Rl^2} g_{\xi \xi \xi} \partial_{\xi}^2\left[(\xi-\tilde{\xi})\Psi(c;\eps)\right] \, d\xi d\tilde{\xi}
\label{defI}
\end{equation}
where
\[
\Psi(c;\eps)=c\phi(c;\eps)-2\Phi(c;\eps).
\]

Expanding the derivatives with respect to $\xi$ in \eq{defI}, using
\eq{defc} to express $c_{\xi \xi}$ in terms of $g_{\xi\xi}$,
and integrating by parts with respect to $\xi$ in the result to remove the third-order derivative
of $g$, we find that $I$ can be expressed as
\[
I = -\frac{5}{2} I_1 + 3 I_2 - I_3
\]
where
\begin{align}
\begin{split}
I_{1}&=\int_{\Rl^2}\Psi^{\prime\prime} (c;\eps) c_{\xi}g^2_{\xi \xi}\, d\xi d\tilde{\xi},
\\
I_{2}&=\int_{\Rl^2} \Psi^{\prime\prime}(c;\eps)c_{\xi}^2 g_{\xi \xi}\, d\xi d\tilde{\xi},
\\
I_{3}&=\int_{\Rl^2} (\xi-\tilde{\xi})\Psi^{\prime\prime\prime}(c;\eps)c_{\xi}^3 g_{\xi \xi}\, d\xi d\tilde{\xi}.
\end{split}
\label{defIj}
\end{align}
The functions $\Psi^{\prime\prime}$, $\Psi^{\prime\prime\prime}$
are given explicitly by
\[
\Psi^{\prime\prime}(c;\eps) = \frac{c}{(1-\eps c)^2},
\qquad
\Psi^{\prime\prime\prime}(c;\eps) = \frac{1 + \eps c}{(1-\eps c)^3}.
\]
In particular, if $|\eps c| \le 1/2$, which is the case if
$g$ satisfies \eq{c2small}, then
\begin{equation}
\label{m}
\vert \Psi^{\prime\prime} (c;\eps)\vert \leq 4 |c|, \qquad \vert \Psi^{\prime\prime\prime}(c;\eps)\vert \leq 12.
\end{equation}
We will estimate the terms in \eq{defIj} separately.

\emph{Estimating $I_{1}$:}
Using (\ref{m}) in \eq{defIj}, we get that
\begin{align}
\begin{split}
\vert I_{1}\vert &\leq 4\int_{\Rl^2} \left|c c_{\xi} g_{\xi\xi}^2\right| \, d\tilde{\xi} d\xi
\\
&\leq 4\sup_{\xi\in \Rl} \left[\int_{\Rl} \left|c c_{\xi}\right| \, d\tilde{\xi}\right]\left(\int_{\Rl}g_{\xi\xi}^2\, d{\xi}\right)
\\
&\leq 4 \sup_{\xi\in \Rl} \left[\left(\int_{\Rl}c^2\, d\tilde{\xi}\right)^{1/2} \left(\int_{\Rl}c_\xi^2\, d\tilde{\xi}\right)^{1/2}\right]
\left(\int_{\Rl}g_{\xi\xi}^2\, d\xi\right).
\end{split}
\label{tempI1}
\end{align}
By a similar argument to the proof of \eq{cl2}, using Taylor's theorem with integral remainder and the Cauchy-Schwartz inequality,
we have from \eq{defc} and \eq{defcxi} that
\begin{align*}
%\|c_\xi\|^2_{L^2} &= \int_{\Rl}c_\xi^2\, d\tilde{\xi}
%\\
\int_{\Rl}c_\xi^2\, d\tilde{\xi} &=\int_{\Rl} \left[\frac{g - \tilde{g} - (\xi-\tilde{\xi}) g_{\xi}}{(\xi-\tilde{\xi})^2}\right]^2\, d\tilde{\xi}
\\
&= \int_0^1 \int_0^1 \int_{\Rl} (1-r)(1-s) g^\prime\left(\tilde{\xi} + r(\xi-\tilde{\xi})\right)
                                           g^\prime\left(\tilde{\xi} + s(\xi-\tilde{\xi})\right)\, d\tilde{\xi} dr ds
\\
&\le \int_0^1 \int_0^1 \int_{\Rl} (1-r)(1-s)
\\
&\qquad
\left(\int_{\Rl} g^{\prime2}\left(\tilde{\xi} + r(\xi-\tilde{\xi})\right)\, d\tilde{\xi}\right)^{1/2}
\left(\int_{\Rl} g^{\prime2}\left(\tilde{\xi} + s(\xi-\tilde{\xi})\right)\, d\tilde{\xi}\right)^{1/2}
dr ds
\\
&\le \left(\int_0^1 \int_0^1 \frac{(1-r)(1-s)}{\sqrt{rs}} \, dr ds\right) \left(\int_{\Rl} g_\xi^2\left(\xi\right)\, d{\xi}\right)
\\
&\le \frac{16}{9} \|g_{\xi\xi}\|^2_{L^2}.
\end{align*}
Thus,
\begin{equation}
\sup_{\xi\in \Rl} \left(\int_{\Rl}c_\xi^2\, d\tilde{\xi}\right)^{1/2} \le \frac{4}{3} \|g_{\xi\xi}\|_{L^2}.
\label{l2cxi}
\end{equation}
Using \eq{cl2} and \eq{l2cxi} in \eq{tempI1}, we get that
\begin{equation*}
\vert I_{1}\vert \leq  A_1 \left\|g_\xi\right\|_{L^2} \|g_{\xi\xi}\|^3_{L^2},
\end{equation*}
where $A_1 = 32/3$ is a numerical constant.

\emph{Estimating $I_{2}$:}
Using (\ref{m}) and \eq{l2cxi} in \eq{defIj}, we get that
\begin{align}
\begin{split}
\vert I_{2}\vert &\leq 4\int_{\Rl^2} \left|c c^2_{\xi} g_{\xi\xi}\right| \, d\tilde{\xi} d\xi
\\
&\leq 4\int_{\Rl} \left(\sup_{\tilde{\xi}\in \Rl}|c|\right) \left(\int_{\Rl} c_{\xi}^2\,d\tilde{\xi}\right)  \left|g_{\xi\xi}\right| \, d{\xi}
\\
&\leq \frac{64}{9} \|g_{\xi\xi}\|^2_{L^2} \int_{\Rl} \left(\sup_{\tilde{\xi}\in \Rl}|c|\right) \left|g_{\xi\xi}\right| \, d{\xi}.
\end{split}
\label{tempI2}
\end{align}
Suppressing the $(t;\eps)$-variables, we observe from \eq{defc} that
\begin{align*}
\sup_{\tilde{\xi}\in \Rl}\vert c\vert &=\sup_{\tilde{\xi}\in \Rl} \left|\frac{g-\tilde{g}}{\xi-\tilde{\xi}}\right|
\\
&= \sup_{\tilde{\xi}\in \Rl}\left|\frac{1}{\xi-\tilde{\xi}}\int^{\xi}_{\tilde{\xi}} g^\prime(z)\, dz \right|
\\
&\le  g^\ast_\xi(\xi),
\end{align*}
where
\[
g^\ast_\xi(\xi) = \sup_{\tilde{\xi}\in \Rl} \frac{1}{|\xi-\tilde{\xi}|} \left|\int_{\tilde{\xi}}^\xi |g^\prime(z)|\, dz\right|
\]
is the maximal function of $g^\prime=g_\xi$, defined using intervals whose left or right endpoint is $\xi$.

Using this inequality and the Cauchy-Schwartz inequality in \eq{tempI2}, we find that
\[
\vert I_{2}\vert
\leq \frac{64}{9} \|g^\ast_\xi\|_{L^2} \|g_{\xi\xi}\|^3_{L^2}.
\]
The maximal operator is bounded on $L^2$, so there exists a numerical constant $M$ such that
\begin{equation}
\|g^\ast_\xi\|_{L^2} \le M\|g_\xi\|_{L^2}.
\label{max_con}
\end{equation}
For example, from \cite{grafakos}, we can take
\[
M = 1 + \sqrt{2}.
\]
It follows that
\[
\vert I_{2}\vert \leq A_2 \|g_\xi\|_{L^2} \|g_{\xi\xi}\|^3_{L^2}
\]
where $A_2 = 64M/9$.

\emph{Estimating $I_{3}$:} Using \eq{defcxi} in \eq{defIj}, we
we can rewrite $I_{3}$ as
\begin{equation*}
I_{3}=\int_{\Rl^2}\Psi^{\prime\prime\prime}(c;\eps) g_{\xi \xi} (c-g_{\xi})c^2_{\xi}\, d\xi d\tilde{\xi}.
\end{equation*}
Splitting this integral into two terms, we get $I_3 = I_3^\prime - I_3^{\prime\prime}$ where
\[
I_{3}^{\prime} =\int_{\Rl^2}\Psi^{\prime\prime\prime}(c;\eps) cc^2_{\xi} g_{\xi \xi}\, d\xi d\tilde{\xi},
\qquad
I_{3}^{\prime\prime} =\int_{\Rl^2}\Psi^{\prime\prime\prime}(c;\eps)c^2_{\xi}g_{\xi}g_{\xi \xi}\, d\xi d\tilde{\xi}.
\]
Using \eq{m}, we have
\[
|I_{3}^{\prime}| \le 12 \int_{\Rl^2} | cc^2_{\xi}g_{\xi \xi}|\, d\xi d\tilde{\xi},
\qquad
|I_{3}^{\prime\prime}| \le 12\int_{\Rl^2} |c^2_{\xi}g_{\xi}g_{\xi \xi}|\, d\xi d\tilde{\xi}.
\]
We estimate $I_{3}^{\prime}$ in exactly the same way as $I_{2}$,
which gives
\[
\vert I_{3}^{\prime}\vert \le A_3^{\prime} \Vert g_{\xi}\Vert _{L^2}\Vert g_{\xi \xi }\Vert^3_{L^2}
\]
where $A_3^{\prime} = 64 M/3$.
We estimate $I_{3}^{\prime\prime}$ in a similar way to $I_{1}$ as
\[
\vert I_{3}^{\prime\prime} \vert
\leq 12 \sup_{\xi\in \Rl}\left[\int_{\Rl} c^2_{\xi}\, d\tilde{\xi}\right]  \left(\int_{\Rl} |g_{\xi}g_{\xi \xi}|\, d\xi\right),
\]
which by use of \eq{l2cxi} and the Cauchy-Schwartz inequality gives
\[
\vert I_{3}^{\prime\prime}\vert \le A_3^{\prime\prime} \Vert g_{\xi}\Vert _{L^2}\Vert g_{\xi \xi }\Vert^3_{L^2}
\]
where $A_3^{\prime\prime} = 64/3$.

Combining these estimates, we get \eq{Iest} with
\begin{equation}
A = 48 + \frac{128}{3} M
\label{Acon}
\end{equation}
where $M$ is the maximal-function constant in \eq{max_con}.
\end{proof}

Using \eq{Iest} in \eq{H2eq}, we find that
\[
\frac{d}{dt} \|g_{\xi\xi}\|_{L^2} \le \frac{1}{2}\eps^2 A  \Vert g_{\xi}\Vert _{L^2}\Vert g_{\xi \xi }\Vert^2_{L^2}.
\]
Since $\|g_\xi\|_{L^2}^2 \le \|g\|_{L^2} \|g_{\xi\xi}\|_{L^2}$ and $\|g\|_{L^2} = \|g_0\|_{L^2}$ is conserved, we get
\begin{equation}
\frac{d}{dt} \|g_{\xi\xi}\|_{L^2} \le \frac{1}{2} \eps^2 A  \Vert g_0\Vert^{1/2}_{L^2}\Vert g_{\xi \xi }\Vert^{5/2}_{L^2}
\label{energyest}
\end{equation}
provided that \eq{c2small} holds. It follows from \eq{energyest} and Gronwall's inequality that if $|\eps|\le \eps_0$, where
$\eps_0$ is sufficiently small, then $\|g_{\xi\xi}\|_{L^2}$ remains finite and \eq{c2small}
holds in some time-interval $0 \le t \le k/\eps^2$, where the constants $\eps_0, k > 0$
may be chosen to depend only on $\|u_0\|_H^2$.
The same estimates hold backward in time, so this completes the proof of Theorem~\ref{th:main}.

By solving the differential inequality \eq{energyest} subject to the constraint \eq{c2small},
we can obtain explicit expressions for $\eps_0$ and $k$. Let
\[
E_0 = \|g_0\|_{L^2}^{1/4} \|g_{0\xi\xi}\|_{L^2}^{3/4},
\]
which is comparable to $\|u_0\|_{H^2}$ from \eq{normest}.
Then we find that Theorem~\ref{th:main} holds with
\[
\eps_0 = \frac{1}{2\sqrt{2} N} \frac{1}{E_0},\qquad
k = \frac{2}{3A} \frac{1}{E_0^2}
\]\
where $N$ is the constant in \eq{gag_nir} and $A$ is the constant in \eq{Acon}.

%------------------------------------------------------------------------------------------------------------------
%------------------------------------------------------------------------------------------------------------------
\section{Normal form transformation}
\label{sec:normal_form}

In this section, we relate the near-identity transformation of the independent variables
used above to a more standard normal form transformation of the dependent variables, of the form
introduced by Shatah \cite{shatah}
\[
v = u + B(u,u)
\]
where $B$ is a bilinear form.

We consider the normal form
transformation $u\mapsto v$ given in \cite{biello}:
\begin{equation}
\label{shortnft}
v=u+\frac{1}{2}\eps \moddf (h^2),\qquad h = \mathbf{H}[u].
\end{equation}
Here, $\df$ denotes the derivative with respect to $x$ and $\moddf = \mathbf{H}\df$.
Differentiating \eq{shortnft} with respect to $t$, using \eq{bheq} to eliminate $u_t$, and simplifying the result, we find that
this transformation removes the nonresonant term of the order $\eps$ from the equation and gives
\begin{equation}
v_{t}+\frac{1}{2}\eps^2 \moddf\,\left[h \moddf(u^2) \right] = \mathbf{H}[v].
\label{shortnfteq}
\end{equation}
The bilinear form $B$ in \eq{shortnft} is not bounded on $H^2$, but
one can show that the normal form transformation \eq{shortnft} is invertible on a bounded set in $H^2$ when $\eps$ is sufficiently small.
We were not able, however, to obtain $H^2$-estimates for $v$ from \eq{shortnfteq}, because
\eq{shortnfteq} contains second-order derivatives, rather than first-order derivatives as in \eq{bheq}, and there is a loss of derivatives
in estimating the $H^s$-norm of $v$. In fact,
for every power of $\eps u$ that one gains
through a normal form transformation of the dependent variable, one introduces an additional derivative.

The appearance of additional derivatives is a consequence of using
a zeroth-order linear term $\mathbf{H}[u]$ to remove a first-order quadratic term $\eps uu_x$.
By contrast, higher-order linear terms lead to normal form transformations that
are easier to analyze. For example, consider the KdV equation
\[
u_t + \eps uu_x =u_{xxx}.
\]
Then, assuming we can ignore difficulties associated with low wavenumbers (\textit{e.g.}\ by considering
spatially periodic solutions with zero mean), we find that the normal form transformation
\[
v = u - \frac{1}{6} \eps \left(\df^{-1} u\right)^2
\]
leads to the equation
\[
v_t - \frac{1}{6}\eps^2 u^2 \left(\df^{-1} u\right)  = v_{xxx}.
\]
In this case, the normal form transformation is bounded and it smooths the nonlinear term.

To explain the connection between the normal form transformation (\ref{shortnft})
and the near-identity transformation \eq{nearid_trans},
we reformulate \eq{shortnft} in a convenient way. Writing $g = \mathbf{H}[v]$ and taking the Hilbert transform
of (\ref{shortnft}), we get the ODE
\begin{equation}
\label{eq1}
g=h-\eps hh_{x}.
\end{equation}
We regard $g(t,x)$ as a given function and use \eq{eq1} to determine the corresponding
function $h$. We may write (\ref{eq1}) as
\begin{equation*}
\frac{h-g}{\eps}-hh_{x}=0,
\end{equation*}
which agrees up to the order $\eps$ with an evolution equation in $\eps$ for $h(t,x;\eps)$:
\begin{equation}
h_{\eps}-hh_{x}=0,\qquad h(t,x;0) = g(t,x).
\label{heps}
\end{equation}
By the method of characteristics, the solution of \eq{heps} is
\[
h(t, x;\eps )=g(t,\xi),\qquad x=\xi-\eps g(t,\xi),
\]
which is the transformation \eq{nearid_trans}. Since \eq{nearid_trans} agrees to the order $\eps$ with a normal form transformation
that removes the order $\eps$ term from \eq{bheq}, this transformation must do so also, as we verified explicitly in Section~\ref{sec:proof}.

It is rather remarkable that the normal form transformation (\ref{shortnft}) can be implemented by making a change of
spatial coordinate in the
equation for $h$, but we do not have a good explanation for why this should be possible.


\begin{thebibliography}{10}
\bibitem{biello} J.~Biello and J.~K.~Hunter, Nonlinear Hamiltonian waves with constant frequency and
surface waves on vorticity discontinuities, \emph{Comm,. Pure Appl. Math.} \textbf{63}, 2009, 303--336.

\bibitem{castro} A.~Castro, D.~C\'ordoba and F.~Gancedo, Singularity formation for a surface wave model,
\emph{Nonlinearity} \textbf{23}, 2010, 2835--2847.

\bibitem{dritschel} D.~G.~Dritschel, The repeated filamentation of two-dimensional vorticity interfaces,
\emph{J. Fluid Mech.} \textbf{194}, 1988,  511--547.

\bibitem{grafakos} L.~Grafakos and S.~Montgomery-Smith,
Best constants for uncentered maximal functions,
\emph{Bull. London Math. Soc.} \textbf{29}, 1997, 60--64.

\bibitem{kato} T.~Kato, The Cauchy problem for quasi-linear symmetric hyperbolic systems, \emph{Arch. Ration. Mech. Anal.} \textbf{58}, 1975, 181--205.

\bibitem{marsden} J.~Marsden and A.~Weinstein, Coadjoint orbits, vortices, and Clebsch variables
for incompressible fluids, \emph{Physica D} \textbf{7}, 1983, 305--323.

\bibitem{shatah} J.~Shatah, Normal forms and quadratic nonlinear Klein-Gordon equations,
\emph{Comm. Pure Appl. Math.} \textbf{38}, 1985, 685--696.

\end{thebibliography}
\end{document}